\documentclass[leqno,twoside]{article} 
\usepackage[latin1]{inputenc}
\usepackage[T1]{fontenc}
\usepackage[english]{babel}
\usepackage{times}
\usepackage{amsmath}
\usepackage[inline]{enumitem}
\usepackage{amssymb} 
\usepackage{graphicx} 
\usepackage{tikz}
\usepackage{bm,dsfont}
\usepackage[margin=1.8cm,top=2.8cm]{geometry}
\usepackage{theorem}
\usepackage{multicol}
\usepackage{float}
\usepackage{hyperref}
\usepackage{cedya2017} % http://cedya2015.uca.es

\newenvironment{demo}{\rm \trivlist \item[\hskip \labelsep{\it
      Proof}.]}{\nopagebreak \hfill $\square$ \endtrivlist}
      
\newtheorem{ejem}[theorem]{Example}

\begin{document}
% >>>>>>>>>>>>>>>>>>>>>>>  Put your title here <<<<<<<<<<<<<<<<<<<<<<<<
\title{%\heading \huge 
\textbf{New uniqueness results for maximal hypersurfaces \\ in Robertson-Walker spacetimes with flat fiber}}

% >>>>>>>>>>>>>>>>>>>>>>> Author's name, affiliation and e-mail <<<<<<<<

 \author{Jos\'e A. S. Pelegr\'in \thanks{Departamento de Geometr\'ia y Topolog\'ia, Universidad de Granada,
Campus Fuentenueva, 18071--Granada(SPAIN).
Email: jpelegrin@ugr.es}}

% >>>>>>>>>>>>>>>>>>>>>>>  Short title and first author for page heading <<<<<<<<<<<<<<<<<<<<
% >>>>>  (notice that if all author's names fit in the line, you may include them all <<<<<<<
\markboth{~\hrulefill\ José A.S. Pelegrín}{Maximal hypersurfaces in Robertson-Walker spacetimes \hrulefill~}

\maketitle

% >>>>>>>>>>>>>>>>>>>>>>>>> Keywords and Abstract <<<<<<<<<<<<<<<<<<<<<
\begin{Abstract}
In this work we study maximal hypersurfaces in spatially open Generalized Robertson-Walker spacetimes with Ricci-flat fiber by means of a generalized maximum principle. In particular, under natural geometric and physical
assumptions we provide new uniqueness and non-existence
results for complete maximal hypersurfaces in spatially open Robertson-Walker spacetimes with flat fiber. Moreover, our results are applied to relevant spacetimes as the steady state spacetime, Einstein-de Sitter spacetime and radiation models.
\vspace*{2ex}\par
\end{Abstract}

% >>>>>>>>>>>>>>>>>>>>>> START OF YOUR PAPER <<<<<<<<<<<<<<<<<<<<<<<<<<<<<<

\begin{multicols}{2}
\section{Introduction}

The importance in General Relativity of maximal and constant mean
curvature spacelike hypersurfaces in spacetimes is well-known; a
summary of several reasons justifying it can be found in \cite{DRT, DRT2, PRR}. Each maximal hypersurface can describe, in some relevant cases, the transition between the expanding and 
contracting phases of a relativistic universe. Moreover, the existence of constant mean curvature (and in particular maximal) hypersurfaces 
is necessary for the study of the structure of singularities in the space of solutions
to the Einstein equations. Also, the deep understanding of this kind of hypersurfaces is essential to prove the positivity of the gravitational mass. 

A maximal hypersurface is (locally) a critical point for 
a natural variational problem, namely of the area functional. From a mathematical point of view, it is necessary to study the maximal hypersurfaces of a spacetime in order to understand its structure. Especially, for some asymptotically flat spacetimes, the existence of a foliation by maximal hypersurfaces is established (see, for instance,
\cite{BF} and references therein). Therefore, their existence and uniqueness appear as crucial issues. Among the uniqueness results, the so-called Calabi-Bernstein theorem stands out. It asserts that the only entire solutions to the maximal hypersurface equation in the Lorentz-Minkowski spacetime, $\mathbb{L}^{n+1}$, are the affine functions defining spacelike hyperplanes \cite{Ca, CY}.

In this paper, basing ourselves on the work done in \cite{PRR2} in order to generalize it, we focus on the problems of uniqueness and non-existence of complete maximal hypersurfaces immersed in a Generalized Robertson-Walker spacetime. The main tool to obtain our results will be the Omori-Yau maximum principle for the Laplacian. In particular, we will be able to deal with maximal hypersurfaces in spatially open Robertson-Walker spacetimes with flat fiber. Note that these models have aroused a great deal of interest, since recent observations have shown that the current universe is very close to a spatially flat geometry \cite{CST}. We will give results that can be used when the fiber is $\mathbb{R}^n$, which is not parabolic for $n \geq 3$ and therefore, cannot be studied in  arbitrary dimension using previous methods \cite{PRR1, RRS}. What is more important, we will not need the hyperbolic angle of the hypersurface to be bounded, which was a restrictive assumption used in previous works studying the spatially open case. Since we are not imposing this restriction, we are able to deal with spacelike hypersurfaces approaching the null Scri boundary at infinity, such as hyperboloids in Minkowski spacetime.
 
Our paper is organized as follows. Section \ref{spre} is devoted to introduce the basic notation used to describe spacelike hypersurfaces in GRW spacetimes. In Section \ref{secsu} we provide an inequality involving the hyperbolic angle of a maximal hypersurface immersed in a GRW spacetime whose fiber is Ricci-flat and obeys the Null Convergence Condition (see Lemma \ref{lemachulo}). This inequality will play a crucial role in our results. Moreover, we also introduce the Omori-Yau maximum principle, which will be key to obtain our uniqueness results. In Section \ref{semr} we get a uniqueness result for complete maximal hypersurfaces (Theorem \ref{cormaxi}). In order to obtain it, the fundamental tool will be a Liouville-type theorem applied to the inequality obtained in  Lemma \ref{lemachulo}. Finally, we give several non-existence results for maximal hypersurfaces in some well-known spacetimes. 

\section{Preliminaries} \label{spre}

Let $(F,g)$ be an $n(\geq 2)$-dimensional (connected) Riemannian manifold, $I$ an open interval in $\mathbb{R}$ and $f$ a positive smooth function defined on $I$. Now, consider the product manifold $\overline{M} = I \times F$ endowed with the Lorentzian metric

\begin{equation}
\label{metr}
\overline{g} = -\pi^*_{_I} (dt^2) +f(\pi_{_I})^2 \, \pi_{_F}^* (g), 
\end{equation}

\noindent where $\pi_{_I}$ and $\pi_{_F}$ denote the projections onto $I$ and
$F$, respectively. The Lorentzian manifold $(\overline{M}, \overline{g})$ is a warped product (in the sense of \cite{O'N}) with base $(I,-dt^2)$, fiber $(F,g)$ and warping function $f$. If we endow $(\overline{M}, \overline{g})$ with the time orientation induced by $\partial_t := \partial / \partial t$ we can call it an $(n+1)$-dimensional Generalized Robertson-Walker (GRW) spacetime. In particular, if the fiber $F$ has constant sectional curvature, it is called a Robertson-Walker spacetime.

The distinguished vector field $K: =~ f({\pi}_I)\,\partial_t$ is timelike and future pointing. From the relation between the
Levi-Civita connection of $\overline{M}$ and those of the base and
the fiber \cite[Cor. 7.35]{O'N}, it follows that

\begin{equation}\label{conexion} 
\overline{\nabla}_X K = f'({\pi}_I)\,X
\end{equation}

\noindent for any $X\in \mathfrak{X}(\overline{M})$, where $\overline{\nabla}$
is the Levi-Civita connection of the Lorentzian metric
(\ref{metr}). Hence, $K$ is conformal and its metrically equivalent $1$-form is closed.

Given an $n$-dimensional manifold $M$, an immersion $\psi: M
\rightarrow \overline{M}$ is said to be spacelike if the
Lorentzian metric (\ref{metr}) induces, via $\psi$, a Riemannian
metric $g_{_M}$ on $M$. In this case, $M$ is called a spacelike
hypersurface. We will denote by $\tau:=\pi_I\circ \psi$ the
restriction of $\pi_I$ along $\psi$.

The time-orientation of $\overline{M}$ allows to take, for each
spacelike hypersurface $M$ in $\overline{M}$, a unique unitary
timelike vector field $N \in \mathfrak{X}^\bot(M)$ globally defined
on $M$ with the same time-orientation as $\partial_t$. Hence, from the wrong way Cauchy-Schwarz inequality, (see \cite[Prop. 5.30]{O'N}) we obtain $\overline{g}(N,\partial_t)\leq -1$ and $\overline{g}(N,\partial_t)= -1$ at a point $p\in M$ if and only if $N(p) = \partial_t(p)$. We will denote by $A$ the shape operator associated to $N$.
Then, the mean curvature function associated to $N$ is given
by $H:= -(1/n) \mathrm{trace}(A)$. As it is well-known, the mean
curvature is constant if and only if the spacelike hypersurface is,
locally, a critical point of the $n$-dimensional area functional for
compactly supported normal variations, under certain constraints of
the volume. When the mean curvature vanishes identically, the
spacelike hypersurface is called a maximal hypersurface.

For a spacelike hypersurface $\psi: M \rightarrow \overline{M}$ with
Gauss map $N$, the hyperbolic angle $\varphi$, at any point
of $M$, between the unit timelike vectors $N$ and $\partial_t$, is
given by $\cosh \varphi=-\overline{g}(N,\partial_t)$. For simplicity,
throughout this paper we will refer to $\varphi$  as the hyperbolic angle function on $M$.

In any GRW spacetime $\overline{M}= I \times_f F$ there
is a remarkable family of spacelike hypersurfaces, namely its
spacelike slices $\{t_{0}\}\times F$, $t_{0}\in I$. It can be easily
seen that a spacelike hypersurface in $\overline{M}$ is a (piece of)
spacelike slice if and only if the function $\tau$ is constant.
Furthermore, a spacelike hypersurface $M$ in $\overline{M}$ is a (piece
of) spacelike slice if and only if $M$ is orthogonal to $\partial_t$. The shape operator of the spacelike slice
$t=t_{0}$ is given by $A=-f'(t_{0})/f(t_{0})\,\mathbb{I}$, where $\mathbb{I}$
denotes the identity transformation, and therefore its (constant)
mean curvature is $H= f'(t_{0})/f(t_{0})$. Thus, a spacelike slice
is maximal if and only if $f'(t_{0})=0$ (and hence, totally
geodesic).

\section{Set up} \label{secsu}

Let $\psi: M \rightarrow \overline{M}$ be an $n$-dimensional
maximal hypersurface immersed in a GRW spacetime
$\overline{M}= I \times_f F$. If we denote by

\[
\partial_t^T:= \partial_t+\overline{g}(N,\partial_t)N
\]
the tangential component of $\partial_t$ along $\psi$, then it is
easy to check that the gradient of $\tau$ on $M$ is

\begin{equation}\label{part}
\nabla \tau=-\partial_t^T
\end{equation}

and so

\begin{equation}\label{sinh}
|\nabla \tau|^2=g_{_M}(\nabla \tau,\nabla \tau)=\sinh^2 \varphi.
\end{equation}

Moreover, since the tangential
component of $K$ along $\psi$ is given by $K^T=K+\overline{g}(K,N)N$, a direct computation from
(\ref{conexion}) gives

\begin{equation}\label{gradcosh}
\nabla \overline{g}(K,N)=-AK^T
\end{equation}

where we have used (\ref{part}), and also

\[
\nabla \cosh \varphi=A\partial_t^T-\frac{f'(\tau)}{f(\tau)}\overline{g}(N,\partial_t)\partial_t^T.
\]

On the other hand, if we represent by $\nabla$ the Levi-Civita
connection of the metric $g_{_M}$, then the Gauss and Weingarten
formulas for the immersion $\psi$ are given, respectively, by

\begin{equation}\label{GF}
\overline{\nabla}_X Y=\nabla_X Y-g_{_M}(AX,Y)N
\end{equation}

and

\begin{equation}\label{WF}
AX=-\overline{\nabla}_X N,
\end{equation}
where $X,Y\in\mathfrak{X}({M})$. Then,  taking the tangential component in
(\ref{conexion}) and using (\ref{GF}) and (\ref{WF}), we get

\begin{equation}\label{KT}
\nabla_X K^T=-f(\tau)\overline{g}(N,\partial_t)AX+f'(\tau)X,
\end{equation}

\noindent where $X\in\mathfrak{X}({M})$ and $f'(\tau):=f'\circ \tau$.

Now we can use (\ref{KT}) as well as the Codazzi equation to obtain that the Laplacian on $M$ of $\cosh \varphi$ is

\begin{gather}
\label{lap5}
\Delta \cosh \varphi = - \overline{\mathrm{Ric}}(\partial_t^T, N) + 2 \frac{f'(\tau)}{f(\tau)} g(A \partial_t^T, \partial_t^T)  \\ \nonumber
 + \cosh \varphi \ \mathrm{trace}(A^2) - \frac{f''(\tau)}{f(\tau)} \cosh \varphi \sinh^2 \varphi \\ \nonumber
+ 3 \frac{f'(\tau)^2}{f(\tau)^2} \cosh \varphi \sinh^2 \varphi + n \frac{f'(\tau)^2}{f(\tau)^2} \cosh \varphi,
\end{gather}

\noindent where $\overline{\mathrm{Ric}}$ represents the Ricci tensor of $\overline{M}$. If we suppose that the fiber $F$ is Ricci-flat and decompose $N$ as $N=N_F-\overline{g}(N,\partial_t)\partial_t$, where
$N_F$ denotes the projection of $N$ on the fiber $F$, we know from \cite[Cor. 7.43]{O'N} that

\begin{equation}
\label{rict}
\overline{\mathrm{Ric}}(\partial_t,\partial_t)=-n \frac{f''(\tau)}{f(\tau)}
\end{equation}

and 

\begin{equation}
\label{ricNF}
\overline{\mathrm{Ric}}(N_F,N_F)= \sinh^2\varphi \left(\frac{f''(\tau)}{f(\tau)}+(n-1)\frac{f'(\tau)^2}{f(\tau)^2} \right).
\end{equation}

Therefore, from (\ref{rict}) and (\ref{ricNF}) we have

\begin{equation}
\label{ritn}
\overline{\mathrm{Ric}}(\partial_t^T, N) = - (n-1) (\log f)''(\tau) \cosh \varphi \sinh^2 \varphi.
\end{equation}

Using (\ref{lap5}), (\ref{ritn}) and the expression of $|\mathrm{Hess}(\tau)|^2$ (see \cite{LR}) we get 

\begin{gather}
\label{clap1}
\cosh \varphi \ \Delta \cosh \varphi = |\mathrm{Hess}(\tau)|^2 \\ \nonumber + n \ \frac{f'(\tau)^2}{f(\tau)^2} \cosh^2 \varphi  
  - \frac{f''(\tau)}{f(\tau)} \cosh^2 \varphi \sinh^2 \varphi \\ \nonumber 
 + 3 \ \frac{f'(\tau)^2}{f(\tau)^2} \cosh^2 \varphi \sinh^2 \varphi 
 - \frac{f'(\tau)^2}{f(\tau)^2} \left( n - 1 + \cosh^ 4 \varphi \right) \\ \nonumber
 - (n-1) (\log f)''(\tau) \cosh^2 \varphi \sinh^2 \varphi.
\end{gather}

Let us assume that our ambient spacetime $\overline{M}$ satisfies the Null Convergence Condition (NCC). It is well known that a spacetime obeys the NCC if and only if its Ricci tensor satisfies $\overline{\mathrm{Ric}}(z,z) \geq 0$ for all lightlike vectors $z$. In particular, in an $(n+1)$-dimensional GRW spacetime with Ricci-flat fiber the NCC is satisfies if and only if $(\log f)'' \leq 0$. Therefore, since $ \frac{1}{2} \Delta \sinh^2 \varphi = \cosh \varphi \ \Delta \cosh \varphi + |\nabla \cosh \varphi|^2$, we can obtain from (\ref{clap1}) the following result

\begin{lemma}
\label{lemachulo}
Let $\psi: M \rightarrow \overline{M}$ be an $n$-dimensional maximal hypersurface immersed in a GRW spacetime
$\overline{M}= I \times_f F$ with Ricci-flat fiber that obeys the Null Convergence Condition, then the hyperbolic angle of the immersion verifies

\begin{equation}
\label{greatresult}
 \frac{1}{2} \Delta \sinh^2 \varphi \geq \left( (n+1)\frac{f'(\tau)^2}{f(\tau)^2} - n \frac{f''(\tau)}{f(\tau)} \right) \sinh^4 \varphi .
\end{equation}

\end{lemma}

\subsection{The Omori-Yau maximum principle}

In order to obtain our results we will make use of the well-known Omori-Yau maximum principle. Following the terminology introduced in \cite{PRS}, we have the following definition.

\begin{definition}
\label{doy}
Let $(M, g)$ be a (not necessarily complete) Riemannian manifold. The Omori-Yau maximum principle for the Laplacian is said to hold on $M$ if for any function $u \in C^2(M)$ with $u^* = \sup u < +\infty$ there exists a sequence of points $\{x_k\}_{k\in \mathbb{N}} \subset M$ satisfying

\begin{center}
\begin{enumerate*}[label=(\roman*),itemjoin={,\quad}]
\item $u(x_k) > u^* - \frac{1}{k}$
\item $|\nabla u (x_k)| < \frac{1}{k}$
\item $\Delta u (x_k) < \frac{1}{k}$
\end{enumerate*}
\end{center}

\noindent for each $k\in \mathbb{N}$. Equivalently, for any function $u \in C^2(M)$ with $u_* = \inf u > -\infty$ there exists a sequence of points $\{x_k\}_{k\in \mathbb{N}} \subset M$ satisfying

\begin{enumerate*}[label=(\roman*),itemjoin={,\quad}]
\item $u(x_k) < u_* + \frac{1}{k}$
\item $|\nabla u (x_k)| < \frac{1}{k}$
\item $\Delta u (x_k) >-\frac{1}{k}$
\end{enumerate*}

\noindent for each $k\in \mathbb{N}$.
\end{definition}

The classical result given by Omori \cite{Om} and Yau \cite{Ya} states that this maximum principle holds on every complete Riemannian manifold with Ricci curvature bounded from below. More generally, it has been shown that a controlled decay of the radial Ricci curvature suffices to guarantee the validity of the Omori-Yau maximum principle on a Riemannian manifold. In particular, the Omori-Yau maximum principle holds on every complete Riemannian manifold whose Ricci curvature has a strong quadratic decay \cite[Thm. 2.2]{CX}, i.e., its Ricci curvature $\mathrm{Ric}$ verifies

\begin{equation}
\label{ricr}
\mathrm{Ric} \geq -c^2 \left( 1 +r^2 \log^2 (2+r) \right),
\end{equation}

\noindent where $c$ is a constant and $r$ is the distance function on the manifold to a fixed point. For maximal hypersurfaces in GRW spacetimes we get the following result

\begin{lemma} \label{ricciacotado1}
Let $\psi: M \rightarrow \overline{M}$ be an $n$-dimensional complete maximal hypersurface immersed in a Robertson-Walker spacetime $\overline{M}= I \times_f F$ with flat fiber that obeys the Null Convergence Condition. Then, the Omori-Yau maximum principle for the Laplacian holds on $M$.
\end{lemma}

\begin{demo}
Given $p\in M$, let us take a local orthonormal frame $\left\{U_1,\ldots,U_n \right\}$ around $p$. From the
Gauss equation we get that the Ricci curvature of $M$, $\mathrm{Ric}$, satisfies

$$\mathrm{Ric}(Y,Y) \geq \sum_k \overline{g}(\overline{\mathrm{R}}(Y,U_k)U_k,Y),
$$ for all $Y \in \mathfrak{X}(M)$, where $\overline{\mathrm{R}}$ denotes the curvature tensor of $\overline{M}$.

Now, from \cite[Prop. 7.42]{O'N} and using the fact that $F$ is flat, we have

\begin{gather}
\label{ricbound}
\sum_k \overline{g}(\overline{\mathrm{R}}(Y,U_k)U_k,Y) = (n-1) \frac{f'(\tau)^2}{f(\tau)^2} |Y|^2 \nonumber \\
 -(n-2)(\log f)''(\tau) \, g(Y,\nabla \tau)^2 - (\log f)''(\tau) |\nabla \tau|^2 |Y|^2 .
\end{gather}

From these equations and our assumptions, we have the Ricci curvature of $M$ to be non-negative. Since $M$ is also complete, the Omori-Yau maximum principle for the Laplacian will hold on $M$.
\end{demo}

Furthermore, if the Omori-Yau maximum principle for the Laplacian holds on a Riemannian manifold, we have the next useful lemma obtained by Nishikawa in \cite{N}.

\begin{lemma}
\label{leman}
Let $M$ be a Riemannian manifold where the Omori-Yau maximum principle for the Laplacian holds and let $u:M\longrightarrow\mathbb{R}$ be a non-negative smooth function on $M$. If there exists a constant $c>0$ such that $\Delta u\geq cu^2$, then $u$ vanishes identically on $M$.
\end{lemma}

\begin{demo}
Let us consider the positive function $F \in C^\infty(M)$ given by $F = \frac{1}{\sqrt{1 + u}}$. If we compute the gradient and the Laplacian of this function we obtain

\begin{equation}
\label{laf}
\Delta u = 6 \frac{|\nabla F|^2}{F^4} - 2 \frac{\Delta F}{F^3}. 
\end{equation}

Since there exists a positive constant $c$ such that $\Delta u\geq cu^2$, we have from (\ref{laf})

\begin{equation}
\label{laf2}
0 \leq c \frac{u^2}{(1+ u^2)^2} \leq 6 |\nabla F|^2 - 2 F \Delta F.
\end{equation}

Applying Omori-Yau maximum principle to the function $F$ we get from (\ref{laf2}) $\sup u = 0$, which finishes the proof.
\end{demo}

\section{Main results} \label{semr}

\begin{theorem}
\label{cormaxi}
Let $\psi: M \rightarrow \overline{M}$ be an $n$-dimensional maximal hypersurface immersed in a GRW spacetime $\overline{M}= I \times_f F$ with Ricci-flat fiber that obeys the Null Convergence Condition. If the Omori-Yau maximum principle for the Laplacian holds on $M$ and $\inf \left\{ (n+1)\frac{f'(\tau)^2}{f(\tau)^2} - n \frac{f''(\tau)}{f(\tau)} \right\}>0$, then $M$ is a spacelike slice $\{t_0\}\times F$ with $f'(t_0)=0$. 
\end{theorem}

\begin{demo}
Using Lemma \ref{leman} in Lemma \ref{lemachulo} we get that $\sinh^2 \varphi$ identically vanishes on $M$.
\end{demo}

As a consequence of Theorem \ref{cormaxi} and Lemma \ref{ricciacotado1} we get

\begin{corollary}
\label{cormaxi2}
Let $\psi: M \rightarrow \overline{M}$ be a complete $n$-dimensional maximal hypersurface immersed in a Robertson-Walker spacetime $\overline{M}= I \times_f F$ with flat fiber that obeys the Null Convergence Condition. If $\inf \left\{ (n+1)\frac{f'(\tau)^2}{f(\tau)^2} - n \frac{f''(\tau)}{f(\tau)} \right\}>0$, then $M$ is a spacelike slice $\{t_0\}\times F$ with $f'(t_0)=0$. 
\end{corollary}

\begin{remark} Observe that the assumption on the function $(n+1)\frac{f'(\tau)^2}{f(\tau)^2} - n \frac{f''(\tau)}{f(\tau)}$ defined on the hypersurface is scarcely restrictive, even if combined with the NCC. In fact, if we consider its extension  $(n+1)\frac{f'(t)^2}{f(t)^2} - n \frac{f''(t)}{f(t)}$ defined on the spacetime, we have from the NCC that $(n+1)\frac{f'(t)^2}{f(t)^2} - n \frac{f''(t)}{f(t)}=\frac{f'(t)^2}{f(t)^2}-n(\log f)''(t)\geq 0$.

However, if we assume that the warping function is defined on the largest possible domain, i.e., it is inextendible, we can find two cases where the required inequality on the infimum does not hold:

\begin{enumerate}
\item When both $f'$ and $f''$ vanish simultaneously at some point in $I=]a, b[$. This obviously happens in the Lorentz-Minkowski spacetime, where an analogous uniqueness result does not hold. Note that $\mathbb{L}^{n+1}$ is a vacuum solution. What is more, if there is real presence of matter in the spacetime we can discard this case.

\item If $\lim\limits_{t \to b} \frac{f'(t)^2}{f(t)^2} = \lim\limits_{t \to b} (\log f)''(t) = 0$. This is the case in the Einstein-de Sitter spacetime. Even more, the inequality will not hold either in the less realistic case where $\lim\limits_{t \to a} \frac{f'(t)^2}{f(t)^2} = \lim\limits_{t \to a} (\log f)''(t) = 0$.
\end{enumerate}

Furthermore, this theorem improves some previous uniqueness results for complete maximal hypersurfaces (see \cite{RRS} and \cite{RRS2}, for instance) without making restrictive assumptions on the maximal hypersurface such as having a bounded hyperbolic angle or lying between two spacelike slices.
\end{remark}

We will give now two models where Theorem \ref{cormaxi} holds.

\begin{ejem}
\label{examp}
\normalfont
Let us consider the Robertson-Walker spacetime $\overline{M} = \mathbb{R} \times_f \mathbb{R}^n$ with warping function $f(t) = e^{-t^2}$.
This spacetime obeys NCC, since $ (\log f)''(t)= -2$.  Moreover, any maximal hypersurface immersed in $\overline{M}$ satisfies

$$\inf \left\{ (n+1)\frac{f'(\tau)^2}{f(\tau)^2} - n \frac{f''(\tau)}{f(\tau)} \right\} = \inf \{ 2n + 4\tau^2 \} > 0.$$
 
\noindent Therefore, the only complete maximal hypersurface in $\overline{M}$ is the spacelike slice $\{0\} \times \mathbb{R}^n$. 

This spacetime models a relativistic universe without singularities (in the sense of \cite[Def. 12.16]{O'N}) that goes from an expanding phase to a contracting one. The physical space in this transition of phase is represented by the spacelike slice $\{0\} \times \mathbb{R}^n$.

\end{ejem}

\begin{ejem}
\label{ext2}
\normalfont
We obtain another example of a Robertson-Walker spacetime satisfying the assumptions in Theorem 4 by considering $\overline{M} = I \times_f \mathbb{R}^n$. Where $I = ]-a, a[$ and the warping function is $f(t) = \sqrt{a^2 - t^2}$, being $a$ a positive constant. Let us remark that this spacetime behaves like the Robertson-Walker model proposed by Friedmann with constant sectional curvature of the fiber equal to one (see \cite[Chap. 12]{O'N}), since it has a big bang singularity at $t=-a$ as well as a big crunch at $t= a$ \cite[Def. 12.16]{O'N}.

For this spacetime, $ (\log f)''(t)= -\frac{a^2+t^2}{(a^2-t^2)^2} \leq 0$, so it satisfies NCC. Furthermore, for every maximal hypersurface in $\overline{M}$ 

$$\inf \left\{ (n+1)\frac{f'(\tau)^2}{f(\tau)^2} - n \frac{f''(\tau)}{f(\tau)} \right\} = \inf \left\{ \frac{n(a^2 + \tau^2) + \tau^2}{(a^2 - \tau^2)^2} \right\} > 0.$$

\noindent Hence, the only complete maximal hypersurface in this spacetime is the spacelike slice $\{0\} \times \mathbb{R}^n$, which represents the physical space in the transition from an expanding phase of the spacetime to a contracting one.

\end{ejem}

Moreover, we can obtain using Corollary \ref{cormaxi2} the following non-existence results for some well-known spacetimes taking into account that a spacelike hypersurface $M$ in a GRW spacetime is said to be bounded away from future infinity when $\sup \tau < +\infty$ on $M$.

\begin{corollary}
\label{coroste}
There are no complete maximal hypersurfaces in the $(n+1)$-dimensional steady state spacetime $\mathbb{R} \times_{e^t} \mathbb{R}^n$.
\end{corollary}

\begin{corollary}
\label{coroein}
There are no complete maximal hypersurfaces  bounded away from future infinity in the $(n+1)$-dimensional Einstein-de Sitter spacetime $\mathbb{R}^+ \times_{t^{2/3}} \mathbb{R}^n$.
\end{corollary}

\begin{corollary}
\label{cororad}
There are no complete maximal hypersurfaces  bounded away from future infinity in the $(n+1)$-dimensional Roberson-Walker Radiation Model spacetime $\mathbb{R}^+ \times_{(2at)^ {1/2}} \mathbb{R}^n$, with $a > 0$.
\end{corollary}

\section*{Acknowledgements}

This research was partially supported by Spanish MINECO and ERDF project MTM2013-47828-C2-1-P.

\label{ultima-pagina:}

\end{multicols}


\begin{thebibliography}{99}
\footnotesize

\bibitem{BF} D. Brill and F. Flaherty,
Isolated maximal surfaces in spacetime, \emph{Commun. Math. Phys.}
\textbf{50} (1984), 157--165.

\bibitem{Ca} E. Calabi, Examples of Bernstein problems for some nonlinear
equations, \emph{P. Symp. Pure Math.}, \textbf{15} (1970),
223--230.

\bibitem{CX} Q. Chen and Y.L. Xin, A generalized maximum principle and its applications in geometry, \emph{Am. J. Math.}, \textbf{114} (1992), 355--366.

\bibitem{CY} S.Y. Cheng and S.T. Yau, Maximal spacelike hypersurfaces in the
Lorentz-Minkowski spaces, \emph{Ann. of Math.},  \textbf{104} (1976),
407--419.

\bibitem{CST} E.J. Copeland, M. Sami and S. Tsujikawa,
Dynamics of dark energy, \emph{Int. J. Mod. Phys. D}, \textbf{15} (2006), 1753--1935.

\bibitem{DRT}  D. de la Fuente, A. Romero and P.J. Torres, Entire spherically spacelike graphs with the prescribed mean curvature function in Schwarzschild and Reissner-Nordstr\"om spacetime, \emph{Classical Quant. Grav.}, \textbf{32} (2015).

\bibitem{DRT2} D. de la Fuente, A. Romero and P.J. Torres, Existence and extendibility of rotationally symmetric graphs with a prescribed higher mean curvature function in Euclidean and Minkowski spaces, \emph{J. Math. Anal. Appl.} \textbf{446} (2017), 1046--1059.

\bibitem{LR} J.M. Latorre and A. Romero, Uniqueness of noncompact spacelike hypersurfaces of constant mean curvature in generalized Robertson-Walker spacetimes, \emph{Geometriae Dedicata}, \textbf{93} (2002), 1--10.

\bibitem{N} S. Nishikawa, On maximal spacelike hypersurfaces in a Lorentzian manifold,
\emph{Nagoya Math. J.}, \textbf{95} (1984), 117--124.

\bibitem{Om} H. Omori, Isometric immersions of Riemannian manifolds, \emph{J. Math. Soc. Japan}, \textbf{19} (1967), 205--214.

\bibitem{O'N} B. O'Neill, \emph{Semi-Riemannian Geometry with
applications to Relativity}, Academic Press, New York, (1983).

\bibitem{PRR} J.A.S. Pelegr\'in, A. Romero and R.M. Rubio, On maximal hypersurfaces in Lorentz manifolds admitting a parallel lightlike vector field, \emph{Classical Quant. Grav.}, \textbf{33} (2016), 055003, 1--8 .

\bibitem{PRR1} J.A.S. Pelegr\'in, A. Romero and R.M. Rubio, On uniqueness of the foliation by comoving observers restspaces of a Generalized Robertson-Walker spacetime, \emph{Gen. Relat. Gravit.}, \textbf{49} (2017), Art. 16, 14pp.

\bibitem{PRR2} J.A.S. Pelegr\'in, A. Romero and R.M. Rubio, Uniqueness of complete maximal hypersurfaces in spatially open $(n+1)$-dimensional Robertson-Walker spacetimes with flat fiber, \emph{Gen. Relat. Gravit.}, \textbf{48} (2016), 1--14.

\bibitem{PRS} S. Pigola, M. Rigoli and A.G. Setti, Maximum principles on Riemannian manifolds and applications, \emph{Mem. Am. Math. Soc.}, \textbf{174} (2005), 99 pp.

\bibitem{RRS} A. Romero, R.M. Rubio and J.J. Salamanca,
Uniqueness of complete maximal hypersurfaces in spatially parabolic
generalized Robertson-Walker spacetimes, \emph{Classical Quant. Grav.}, \textbf{30} (2013) 115007--115020.

\bibitem{RRS2} A. Romero, R.M. Rubio and J.J. Salamanca, Complete maximal hypersurfaces in certain spatially open generalized Robertson-Walker spacetimes, \emph{RACSAM Rev. R. Acad. A.}, \textbf{109} (2015), 451--460.

\bibitem{Ya} S.T. Yau, Harmonic functions on complete Riemannian manifolds, \emph{Comm. Pure Appl. Math.}, \textbf{28} (1975), 201-228.


\end{thebibliography}
\end{document}